\documentclass[12pt]{article}
\usepackage{amssymb}
\usepackage{amsmath}
\usepackage{graphicx}
\usepackage{color}

\newtheorem{thm}{Theorem}[section]
 \newtheorem{cor}[thm]{Corollary}
 \newtheorem{lem}[thm]{Lemma}
 \newtheorem{prop}[thm]{Proposition}
 \newtheorem{defn}[thm]{Definition}
 \newtheorem{rem}[thm]{Remark}
 
 \numberwithin{equation}{section}

 \newcommand{\rbx}{\hfill{\rule{1ex}{1ex}}}

\newcommand{\ind}{\mathrm{ind}\,}

\newcommand{\diag}{\mathrm{diag}\,}

 \newcommand{\vp}{\varphi}

\newcommand{\coker}{\mbox{\rm coker}\,}

\newcommand{\im}{\mathrm{im}\,}

\newcommand{\cL}{\mathcal{L}}

\newcommand{\cU}{\mathcal{U}}

\newcommand{\sC}{{\mathbb C}}

\newcommand{\sR}{{\mathbb R}}

\newcommand{\sT}{{\mathbb T}}

\begin{document}

\vspace*{10mm}

\begin{center}
{\Large\textbf{Some classes of Wiener--Hopf plus Hankel operators
and the Coburn-Simonenko Theorem}}\footnote{This research was
supported by the Universiti Brunei Darussalam under Grant
UBD/GSR/S\&T/19.}
\end{center}

\vspace{5mm}

%\author{Victor D. Didenko}
\begin{center}

\textbf{Victor D. Didenko and Bernd Silbermann}

 %\title[\textit{Invertibility of Toeplitz plus Hankel operators}]{}

\vspace{2mm}


%    author one information
 Universiti Brunei Darussalam,
Bandar Seri Begawan, BE1410  Brunei; diviol@gmail.com

Technische Universit{\"a}t Chemnitz, Fakult{\"a}t f\"ur Mathematik,
09107 Chemnitz, Germany; silbermn@mathematik.tu-chemnitz.de

 \end{center}

%\keywords{Toeplitz plus Hankel operator, Invertibility}

  \vspace{10mm}

\textbf{2010 Mathematics Subject Classification:} Primary 47B35,
47B38; Secondary 47B33, 45E10

\textbf{Key Words:} Wiener--Hopf plus Hankel operator,
Coburn--Simonenko theorem, invertibility

%\dedicatory{Communicated by }
%% We use this for communication information.
%% If you want do dedicate your paper to somebody, then please use \thanks{}

%\begin{document}

%{\begin{flushleft}\baselineskip9pt\scriptsize
%PUBLICATIONS DE L'INSTITUT MATH\'EMATIQUE\newline
%Nouvelle s\'erie, tome 87(101) (2010), od--do \hfill DOI:
%\end{flushleft}}
\vspace{18mm} \setcounter{page}{1} \thispagestyle{empty}

\begin{abstract}
Wiener--Hopf plus Hankel operators $W(a)+H(b):L^p(\sR^+)\to
L^p(\sR^+)$ with generating functions  $a$ and $b$ from a subalgebra
of $L^\infty(\sR)$ containing almost periodic functions and Fourier
images of $L^1(\sR)$-functions are studied. For $a$ and $b$
satisfying the so-called matching condition
 \begin{equation*}
a(t) a(-t)=b(t) b(-t),\quad t\in\sR,
 \end{equation*}
we single out some classes of operators $W(a)+H(b)$ which are
subject to Coburn--Simonenko theorem.
\end{abstract}

 \sloppy
\section{Introduction\label{s1}}

The classical Coburn--Simonenko Theorem states that for a Toeplitz
or Wiener--Hopf operator $A$ with a scalar nonzero generating
function, at least one of the numbers $\dim\ker A$ or $\dim\coker A$
is equal to zero. Thus if it is known that the corresponding
operator is Fredholm with index zero, the Coburn--Simonenko Theorem
implies that this operator is invertible. Note that Fredholmness of
such operators with generating functions from various classes is
well understood. On the other hand, for Toeplitz plus Hankel
operators $T(a)+H(b)$ with piecewise continuous generating functions
$a$ and $b$ their Fredholm properties can be derived by a direct
application of results \cite[Sections 4.95--4.102]{BS:2006},
\cite[Sections 4.5 and 5.7]{RSS:2011}, \cite{RS1990}. The case of
quasi piecewise continuous generating functions has been studied in
\cite{Si:1987}, whereas formulas for the index of the operators
$T(a)+H(b)$ considered on various Banach and Hilbert spaces and with
various assumptions about the generating functions $a$ and $b$ have
been established in \cite{DS:2013,RS:2012}. It is also worth
mentioning that lately a lot of effort has been spent to obtain
information concerning the kernel and cokernel dimensions of
Toeplitz plus Hankel or Wiener--Hopf plus Hankel operators. Here we
are not going to discuss the history of these investigations in much
detail, but just mention a few important developments. For example,
in the works of Ehrhardt \cite{Ehr:2004h, E2004} and Ehrhardt and
Basor \cite{BE2004, BE2006,BE:2013}, Toeplitz plus Hankel operators
have been studied in $H^p$-spaces on the unit circle $\sT$ mainly
under the assumption that the generating functions of these
operators are piecewise continuous, satisfy an algebraic relation,
and that the operators are Fredholm. Wiener--Hopf plus Hankel
operators have received less attention in the literature and results
are scarce (see, for example, \cite{CS:2012} and references there).
In addition, in most cases the conditions imposed on the generating
functions are very restrictive and ensure that the problem can be
handled in a more or less straightforward way.

Let us now describe the problem studied in the present paper.
Consider the set $G$ of all functions of the form
\begin{equation}\label{cst1}
a(t)=\sum_{j=-\infty}^\infty a_j e^{i\delta_j
t}+\int_{-\infty}^\infty k(s)e^{its}\,ds, \quad-\infty<t<\infty,
\end{equation}
where $\delta_j\in\sR$ are pairwise distinct and
\[
\sum_{j=-\infty}^\infty |a_j|<\infty, \quad \int_{-\infty}^\infty
|k(s)|\,ds<\infty.
\]
The set $G$ actually forms a commutative unital Banach algebra under
pointwise operations and the norm
\[
\|a\|:=\sum_{j=-\infty}^\infty
|a_j|+\int_{-\infty}^\infty|k(s)|\,ds.
\]
This algebra $G$ contains both the algebra $AP_w$ of all almost
periodic functions with absolutely convergent Fourier series and the
algebra $\cL_0$ of all Fourier transforms of functions from
$L^1(\sR)$. Moreover, the algebra $G$ is the direct sum of $AP_w$
and $\cL_0$, and $\cL_0$ is an ideal in $G$. A function $a\in G$ is
invertible in $G$ if and only if it satisfies the condition
$\inf_{t\in\sR} |a(t)|>0$. Moreover, if $b\in AP_w$, $k\in \cL_0$,
and $b+k$ is invertible in $G$, then $b$ is also invertible in
$AP_w$ (see \cite[Chapter VII]{GF1974}). Further, let us introduce
the subalgebra $G^+$ ($G^-$) of the algebra $G$, which consists of
all functions \eqref{cst1} such that all numbers $\delta_j$ are
nonnegative (nonpositive) and function $k$ vanishes on the negative
(positive) semi-axis. It is clear that the functions from $G^+$ and
$G^-$ admit holomorphic extensions to the upper and to the lower
half-plane, respectively, and the intersection of the sets $G^+$ and
$G^-$ contains constant functions only.

If $b\in AP_w$, $k\in\cL_0$, and the element $a=b+k$ is invertible
in $G$, then the numbers
\[
\nu(a):=\lim_{l\to\infty}\frac{1}{2l} [\arg b(t)]_{-l}^l,\quad
\text{and}\quad n(a):=\frac{1}{2\pi} [\arg
(1+b^{-1}(t)k(t)]_{t=-\infty}^\infty,
\]
are well defined. In particular, the first limit exists because $b$
is an almost periodic function.

Let $\sR^+:=(0,\infty)$ and let $P$ be the projection operator from
$L^p(\sR)$, $1\leq p\leq \infty$ onto $L^p(\sR^+)$, that is
$P:f\mapsto f|_{\sR^+}$. Analogously, $Q$ is the projection operator
from $L^p(\sR)$ onto $L^p(\sR^-)$, $\sR^-:=(-\infty,0)$. In what
follows we will identify the space $L^p(\sR^+)$ ($L^p(\sR^-)$) with
the subspace of $L^p(\sR)$ consisting of all functions vanishing on
$\sR^-$ ($\sR^+$). Note that $P^2=P$ and $Q^2=Q$.

Each function $a\in G$,
\[
a(t)=\sum_{j=-\infty}^\infty a_j e^{i\delta_j
t}+\int_{-\infty}^\infty k(s) e^{its} \, ds,
\]
generates two operators $W^0(a):L^p(\sR)\to L^p(\sR)$ and
$W(a):L^p(\sR^+)\to L^p(\sR^+)$ defined by
\begin{align}\label{cst2}
(W^0(a)f)(t)&:=\sum_{j=-\infty}^\infty a_j f(t-\delta_j)+\int_{-\infty}^\infty k(t-s) f(s)\,ds, \\
W(a)f&:=PW^0(a) f.\notag
\end{align}
These operators belong to the spaces $\cL (L^p(\sR))$ and
$\cL(L^p(\sR^+))$, respectively, i.e., they are linear bounded
operators. Moreover, the mappings $G\to\cL (L^p(\sR))$ and
$G\to\cL(L^p(\sR^+))$ defined, respectively, by
\[
a\mapsto W^0(a)\, \text{ and }\, a\mapsto W(a) \,,
\]
are injective linear bounded mappings. The function $a$ is referred
to as the generating function, or the symbol, for both operators
$W^0(a)$ and $W(a)$. The Fredholm theory for the operators $W^0(a)$,
$a\in G$ is relatively simple. An operator $W^0(a)$ is semi-Fredholm
if and only if $a$ is invertible in $G$. A proof of this result is
implicitly contained in the proof of Theorem 2.4, \S 2, Chapter VII
in \cite{GF1974}.

Note that the convolution operator \eqref{cst2} is shift invariant
that is $W^0(a)\tau_v=\tau_v W^0(a)$ for any $v\in\sR$, where
$\tau_v$ is the operator defined by $(\tau_v f)(t):=f(t-v)$. The
operator $W(a)$ is called integro-difference operator
\cite[Chapter~VII]{GF1974}. It is shown in \cite[Sections 9.4 and
9.21]{BS:2006} that integro-difference operators are indeed
Wiener--Hopf integral operators. If $a$ does not vanish identically,
then $W(a)$ has a trivial kernel or a dense range in $L^p(\sR^+)$ at
least for $1<p<\infty$ and this is the Coburn--Simonenko Theorem for
such class of operators (see \cite[Section 9.5 (d)]{BS:2006}).

Now we can formulate the following result.

\begin{thm}[Gohberg/Feldman \cite{GF1974}]\label{t1}
If $a\in G$, then the operator $W(a)$ is one-sided invertible in
$L^p(\sR^+)$ for $1\leq p\leq \infty$ if and only if $a$ is
invertible in $G$. Further, if $a\in G$ is invertible in $G$, then
the following assertions are true.

\begin{enumerate}
\item If $\nu(a)>0$, then the operator $W(a)$ is invertible from the left and $\dim\coker W(a)=\infty$.
\item If $\nu(a)<0$, then the operator $W(a)$ is invertible from the right and $\dim\ker W(a)=\infty$.
\item If $\nu(a)=0$, then the operator $W(a)$ is invertible from the left (right) if $n(a)\geq 0$ ($n(a)\leq 0$) and
\[
\dim \coker W(a)=n(a) \quad (\dim \ker W(a)=-n(a)).
\]
\item If $a\in G$ is not invertible in $G$, then $W(a)$ is not a semi-Fredholm operator.
\end{enumerate}
\end{thm}

\begin{rem}\label{r1}
Using the Coburn--Simonenko Theorem, one can show that if $W(a)$ is
normally solvable and $a\neq 0$, then $a$ is invertible in $G$, at
least in the case where the operator $W(a)$ acts on the space
$L^p(\sR^+)$, $p\in (1,\infty)$.
\end{rem}
Let us introduce Hankel operators. For, consider the operator
$J:L^p(\sR)\to L^p(\sR)$ defined by $J\vp :=\widetilde{\vp}$, where
$\widetilde{\vp}(t):=\vp(-t)$. If $a\in G$ and $1\leq p\leq\infty$,
then on the space $L^p(\sR^+)$ the Hankel operators $H(a)$ and
$H(\widetilde{a})$ are defined as follows
\begin{align*}
H(a)&:\vp\mapsto PW^0(a)QJ\vp, \\
H(\widetilde{a})&: \vp\mapsto JQW^0(a)P\vp.
\end{align*}
Note that $JQ W^0(a)P=PW^0(\widetilde{a})QJ$, and the last identity
is the consequence of the following relations
\begin{equation}\label{cst3}
J^2=I,\quad JQ=PJ, \quad JP=QJ,\quad JW^0(a)J=W^0(\widetilde{a}).
\end{equation}

On the space $L^p(\sR)$, $1\leq p\leq\infty$ we also consider the
operators $\cU$ and $\cU^{-1}$ defined by
\begin{align*}
(\cU\vp)(t)&:=\vp(t)-2\int_{-\infty}^t e^{s-t}\vp(s)\,ds, \quad-\infty<t<\infty, \\
(\cU^{-1}\vp)(t)&:=\vp(t)-2\int_t^{\infty} e^{t-s}\vp(s)\,ds,
\quad-\infty<t<\infty.
\end{align*}
It is well known \cite{GF1974} that
\[
\cU=W^0(\chi), \quad \cU^{-1}=W^0(\chi^{-1}),
\]
where $\chi(t):=(t-i)/(t+i)$, $\chi^{-1}(t):=(t+i)/(t-i)$,
$t\in\sR$. Moreover, since
$W^0(\chi)W^0(\chi^{-1})=W^0(\chi\chi^{-1})$, we get $\cU
\cU^{-1}=\cU^{-1}\cU=I$.

One of the aims of this work is to establish a Coburn--Simonenko
Theorem for the operators $W(a)+H(a\chi)$ and $W(a)-H(a\chi^{-1})$,
where $a\in G$ is invertible. Recall that the semi-Fredholmness of
the operators $W(b)+H(c):L^p(\sR^+)\to L^p(\sR^+)$, $b,c\in G$
implies that the element $b$ is invertible in $G$ at least in the
case where $1<p<\infty$. Indeed, the proof of Theorem 2.30 in
\cite{BS:2006} with the shifts $U^{\pm n}$ and the Toeplitz
operators $T(a)$ replaced, respectively, by the translations
$\tau_{\pm \nu}$, $\nu\in \sR^+$ and the operators $W(b)+H(c)$
implies that $\|W^0(b)f\| \geq c\|f\|\;\text{ for all } f\in
L^p(\sR)$. But then $W^0(b)$ is semi-Fredholm and, therefore, $b$ is
invertible in $G$. For $p=1$ this proof does not work. Nevertheless,
we conjecture that for $p=1$ the result is also true. Therefore, the
above requirement of the invertibility of the element $a$ is not too
restrictive.

Finally, let us also mention that if $a,b\in G$, then $W^0(a
b)=W^0(a) W^0(b)$, and if $a\in G^-, c\in G^+$, and $b\in G$, then
$W(abc)=W(a) W(b) W(c)$. Moreover, in the following we will make use
of the identities
\begin{equation}\label{cst4}
\begin{aligned}
W(ab)&=W(a)W(b)+H(a) H(\widetilde{b}),\\
H(ab)&=W(a)H(b)+H(a)W(\widetilde{b}).
\end{aligned}
\end{equation}

\section{Kernels of Wiener--Hopf plus Hankel operators.\\General properties\label{s2}}

In this section we establish certain relations between the kernels
of Wiener--Hopf plus Hankel operators and matrix Wiener-Hopf
operators in the case where the generating functions $a,b\in G$. The
corresponding results for Toeplitz plus Hankel operators
$T(a)+H(b)$, $a,b\in L^\infty$ have been obtained recently
\cite{DS:2013e}. Taking into account Theorem \ref{t1} we can always
assume that $a$ is invertible in~$G$. Along with the operator
$W(a)+H(b)$ let us also consider the Wiener--Hopf minus Hankel
operator $W(a)-H(b)$ and the Wiener--Hopf operator $W(V(a,b)))$
defined by the matrix
\[
V(a,b):=
\begin{pmatrix}
a-b\widetilde{b}\widetilde{a}^{-1} & d \\
-c & \widetilde{a}^{-1}\\
\end{pmatrix},
\]
where $c:=\widetilde{b}\widetilde{a}^{-1}$,
$d:=b\widetilde{a}^{-1}$.

The following lemma describes connections between the solutions of
homogeneous equations with Wiener--Hopf plus/minus Hankel operators
and the solutions of the associated homogeneous equation with the
matrix Wiener--Hopf operator $W(V(a,b))$.

\begin{lem}\label{l1}
Assume that $a,b\in G$, $a$ is invertible in $G$, and the operators
$W(a)\pm H(b)$ are considered on the space $L^p(\sR^+)$, $1\leq
p\leq\infty$.
\begin{itemize}
\item If $(\vp,\psi)^T\in \ker W(V(a,b))$, then
\begin{multline}\label{cst5}
\!\!\!\!\!\!\!\!(\Phi,\!\! \Psi)^T\!\!=\!\!\frac12\big(\vp-J Q
W^0(c)\vp+J QW^0(\widetilde{a}^{-1})\psi,\vp+J Q W^0(c)\vp
-J QW^0(\widetilde{a}^{-1})\psi\big)^T\\
\in\ker\diag\big(W(a)+H(b),W(a)-H(b)\big)
\end{multline}
\item If $(\Phi,\Psi)^T\in\ker\diag\big(W(a)+H(b),W(a)-H(b)\big)$, then
\begin{equation}\label{cst6}
\big(\Phi+\Psi,P(W^0(\widetilde{b})(\Phi+\Psi)+W^0(\widetilde{a})JP(\Phi-\Psi))\big)^T\in\ker
W(V(a,b)).
\end{equation}
\end{itemize}
Moreover, the operators
\begin{align*}
E_1:\ker W(V(a,b)) \to \ker\diag\big(W(a)+H(b),W(a)-H(b)\big),\\%[1ex]
E_2:\ker\diag\big(W(a)+H(b),W(a)-H(b)\big) \to \ker W(V(a,b)),
\end{align*}
defined, respectively, by relations \eqref{cst5} and \eqref{cst6}
are mutually inverse.
\end{lem}

\textbf{Proof.} Consider the operators
\begin{align}\label{cst7}
A&:=
\begin{pmatrix}I & 0\\W^0(\widetilde{b}) & W^0(\widetilde{a})\end{pmatrix}
\begin{pmatrix}I & I\\J &-J\end{pmatrix}, \qquad
B_1 :=2\begin{pmatrix}I & J \\I &-J \end{pmatrix},\\
B_2&:=\diag(I,I)-\diag(P,Q)
\begin{pmatrix}W^0(a) & W^0(b)\\W^0(\widetilde{b}) & W^0( \widetilde{a})\end{pmatrix}\diag(Q,P), \notag \\
B_3 &:=\diag(I,I)+\diag(P,P)
\begin{pmatrix}
W^0(a-b\widetilde{b}\widetilde{a}^{-1}) & W^0(d)\\-W^0(c) &
W^0(\widetilde{a}^{-1})\end{pmatrix}\diag(Q,Q).\notag
\end{align}
Elementary but tedious computations show that the the operator
\[
\diag\big(W(a)+H(b)+Q, W(a)-H(b)+Q\big)
\]
can be represented as the product of three matrix operators, viz.
\begin{equation}\label{cst8}
\begin{pmatrix}
W(a)+H(b)+Q & 0\\0 & W(a)-H(b)+Q\end{pmatrix}=B(W(V(a,b)))+\diag
(Q,Q)A\,,
\end{equation}
where $B:=B_1 B_2 B_3$. The operator $A:L^p(\sR)\times L^p(\sR)\to
L^p(\sR)\times L^p(\sR)$ is invertible because $a$ is invertible in
$G$, and it is well known that all the operators $B_1,B_2,B_3$ are
invertible as well. Therefore, relations \eqref{cst7}--\eqref{cst8}
imply that for any $(\vp,\psi)^T\in\ker W(V(a,b)))$, the element
$A^{-1}((\vp,\psi)^T)$ belongs to the set
\begin{align*}%\label{eq}
&\ker\diag\big(W(a)+H(b)+Q,W(a)-H(b)+Q\big)\\
 &\qquad\quad =\ker\diag\big(W(a)+H(b),W(a)-H(b)\big)
\end{align*}.
Hence
\[
\diag(P,P)A^{-1}((\vp,\psi)^T)=A^{-1}((\vp,\psi)^T).
\]
Computing the left-hand side of the last equation, one obtains
relation \eqref{cst5}. Analogously, if
$(\Phi,\Psi)^T\in\ker\diag\big(W(a)+H(b),W(a)-H(b)\big)$, then
$A((\Phi,\Psi)^T)\in\ker W(V(a,b))$ and
$\diag(P,P)\,A((\Phi,\Psi)^T)=A((\Phi,\Psi)^T)$, so representation
\eqref{cst6} follows.

Now let $(\vp,\psi)$ and $(\Phi,\Psi)$ be as above. Then
\begin{align*}
\diag (P,P)A\diag(P,P)A^{-1}((\vp,\psi)^T)&=AA^{-1}((\vp,\psi)^T),\\
\diag(P,P)A^{-1}\diag(P,P)A((\Phi,\Psi)^T)&=A^{-1}A((\Phi,\Psi)^T),
\end{align*}
which completes the proof. \rbx

From now on we will always assume that the generating functions $a$
and $b$ satisfy the condition
\begin{equation}\label{cst9}
a \widetilde{a}=b \widetilde{b}.
\end{equation}
Analogously to \cite{DS:2013d}, relation~\eqref{cst9} is called
matching condition, and if $a$ and $b$ satisfy \eqref{cst9}, then
the duo $(a,b)$ is called a matching pair. For each matching pair
$(a,b)$ one can assign another matching pair $(c,d)$ with
$c:=\widetilde{b}\widetilde{a}^{-1}$ and $d:=b \widetilde{a}^{-1}$.
Such a pair $(c,d)$ is called the subordinated pair for $(a,b)$, and
it is easily seen that the functions which constitutes a
subordinated pair have a specific property, namely
$c\tilde{c}=1=d\widetilde{d}$. Throughout this paper any function
$g\in G$ satisfying the condition
\[
g \widetilde{g}=1,
\]
is called matching function. In passing note that the matching
functions $c$ and $d$ can also be expressed in the form
\[
c=ab^{-1}, \quad d=\widetilde{b}^{-1} a.
\]
Besides, if $(c,d)$ is the subordinated pair for a matching pair
$(a,b)$, then $(\overline{d},\overline{c})$ is the subordinated pair
for the matching pair $(\overline{a},\overline{\widetilde{b}})$
which defines the adjoint operator
\begin{equation}\label{cst10}
(W(a)+H(b))^*=W(\overline{a})+H(\overline{\widetilde{b}})
\end{equation}
for the operator $W(a)+H(b)$. Further, a matching pair $(a,b)$ is
called Fredholm, if the Wiener--Hopf operators $W(c)$ and $W(d)$ are
Fredholm.

If $(a,b)$ is a matching pair, then the corresponding
matrix--function $V(a,b)$ takes the form
\[
V(a,b)=
\begin{pmatrix}0 & d \\-c & \widetilde{a}^{-1}\end{pmatrix},
\]
where $(c,d)$ is the subordinated pair for the pair $(a,b)$.
Moreover, similarly to the corresponding representation of the
matrix Toeplitz operator $T(V(a,b))$ from \cite{DS:2013d}, the
operator $W(V(a,b)))$ can be represented as the product of three
matrix Wiener--Hopf operators
\begin{align} \label{cst11}
W(V(a,b)))&=\begin{pmatrix}0 & W(d) \\-W(c) &
W(\widetilde{a}^{-1})\end{pmatrix}
\\
&=\begin{pmatrix}-W(d) & 0 \\0 & I\end{pmatrix}
\begin{pmatrix}0 &-I \\I & W(\widetilde{a}^{-1})\end{pmatrix}
\begin{pmatrix}-W(c) & 0\\0 & I\end{pmatrix},\notag
\end{align}
where the operator
\[
D:=\begin{pmatrix}0 &-I \\I & W(\widetilde{a}^{-1})\end{pmatrix}
\]
in the right-hand side of \eqref{cst11} is invertible and
\[
D^{-1}=\begin{pmatrix}W(\widetilde{a}^{-1}) & I \\-I &
0\end{pmatrix}.
\]
Note that a useful representation for the kernel of the block
Toeplitz operator $T(V(a,b)))$ defined by a matching pair $(a,b)$,
has been derived recently. Following
\cite[Proposition~3.3]{DS:2013e}, one can also obtain a similar
result for the block Wiener-Hopf operator $W(V(a,b))$.

\begin{prop}\label{p3.4}
Let $(a,b)\in G\times G$ be a matching pair such that the operator
$W(c)$, $c=ab^{-1}$, is invertible from the right. Then
\[
\ker W(V(a,b)))=\Omega(c)\dotplus\widehat{\Omega}(d)
\]
where
\begin{align*}
\Omega(c)&:=\big\{(\vp,0)^T:\vp\in\ker W(c)\big\},\\
\widehat{\Omega}(d)&:=\big\{(W_r^{-1}(c)W(\widetilde{a}^{-1})s,s)^T:s\in\ker
W(d)\big\},
\end{align*}
and $W_r^{-1}(c)$ is one of the right inverses for the operator
$W(c)$.
\end{prop}

\textbf{Proof.} It is clear that $\Omega(c)$ and
$\widehat{\Omega}(d)$ are closed subspaces of $\ker W(V(a,b))$ and
$\Omega(c) \cap \widehat{\Omega}(d)=\{0\}$.

If $(y_1,y_2)^T\in\ker W(V(a,b))$, then $W(d)y_2=0$, and
$W(c)y_1=W(\widetilde{a}^{-1})y_2$. Since $W_r^{-1}(c)$ is
left-invertible, the space $L^p(\sR^+)$ is the direct sum of the
closed subspaces $\ker W(c)$ and $\im W_r^{-1}(c)$, i.e.,
$L^p(\sR^+)=\ker W(c)\dotplus \im W_r^{-1}(c)$. Consequently, the
element $y_1$ can be represented in the form $y_1=y_{10}+y_{11}$,
where $y_{10}\in\ker W(c)$ and $y_{11}\in\im W_r^{-1}(c)$. Moreover,
there is a unique vector $y_3\in L^p(\sR^+)$ such that
$y_{11}=W_r^{-1}(c)y_3$, so we get
\[
W(c) y_1=W(c)(W_r^{-1}(c)y_3+y_{10})=y_3=W(\widetilde{a}^{-1})y_2.
\]
It implies that $y_1=W_r^{-1}(c)W(\widetilde{a}^{-1})y_2+y_{10}$,
what leads to the representation
\[
(y_1,y_2)^T=(W_r^{-1}(c)W(\widetilde{a}^{-1})y_2,y_2)^T+(y_{10},0)^T,
\]
with $(W_r^{-1}(c)
W(\widetilde{a}^{-1})y_2,y_2)^T\in\widehat{\Omega}(d)$ and
$(y_{10},0)^T\in \Omega(c)$. \rbx

Thus $\vp \in \ker W(c)$ implies that $(\vp,0)^T\in\ker W(V(a,b)))$
and by Lemma \ref{l1}
\begin{equation}\label{cst12}
\begin{aligned}
\vp-JQW^0(c)P\vp\in\ker(W(a)+H(b)),\\
\vp+JQW^0(c)P\vp\in\ker(W(a)-H(b)).
\end{aligned}
\end{equation}
It is even more remarkable that the functions $\vp-JQW^0(c)P\vp$ and
$\vp+JQW^0(c)P\vp$ belong to the kernel of the operator $W(c)$ as
well.

\begin{prop}\label{p1}
Let $g\in G$ be a matching function, i.e., $g \widetilde{g}=1$. If
$f\in\ker W(g)$, then $JQW^0(g)Pf\in\ker W(g)$ and $(J
QW^0(g)P)^2f=f$.
\end{prop}

\textbf{Proof.} If $g \widetilde{g}=1$ and $f\in \ker W(g)$, then
\begin{align*}
W(g)(JQW^0(g)Pf)&=PW^0(g)PJ QW^0(g)Pf=JQW^0(\widetilde{g})QW^0(g)Pf\\
&=JQW^0(\widetilde{g})W^0(g)Pf-JQW^0(\widetilde{g})PW^0(g)Pf=0,
\end{align*}
and assertion (i) follows. On the other hand, for any $f\in\ker
W(g)$ one has
\begin{align*}
(JQW^0(g)P)^2f&=JQW^0(g)PJQW^0(g)Pf=PW^0(\widetilde{g})QW^0(g)Pf&\\
&=PW^0(\widetilde{g})W^0(g)Pf-PW^0(\widetilde{g})PW^0(g)Pf=f,
\end{align*}
which completes the proof. \rbx

Consider now the operator $\mathbf{P}(g):=JQW^0(g)P\big|_{\ker
W(g)}$. Proposition~\ref{p1} implies that $\mathbf{P}(g):\ker
W(g)\to\ker W(g)$ and $\mathbf{P}^2(g)=I$. Thus on the space $\ker
W(g)$ the operators $\mathbf{P}^-(g):=(1/2)(I-\mathbf{P}(g))$ and
$\mathbf{P}^+(g):=(1/2)(I+\mathbf{P}(g))$ are complementary
projections generating a decomposition of $\ker W(g)$. Moreover,
relations \eqref{cst12} lead to the following result.

\begin{cor}\label{c1}
Let $(c,d)$ be the subordinated pair for a matching pair $(a,b)\in
G\times G$. Then $\ker
W(c)=\im\mathbf{P}^-(c)\dotplus\im\mathbf{P}^+(c)$, and the
following relations hold
\begin{equation}\label{cst13}
\im\mathbf{P}^-(c)\subset\ker(W(a)+H(b)) ,\quad
\im\mathbf{P}^+(c)\subset\ker(W(a)-H(b)).
\end{equation}
\end{cor}

Relations \eqref{cst13} show the influence of the operator $W(c)$ on
the kernels of the operators $W(a)+H(b)$ and $W(a)-H(b)$. Let us now
clarify the role of another operator--viz.~the operator $W(d)$, in
the structure of the kernels of the operators $W(a)+H(b)$ and
$W(a)-H(b)$. Assume additionally that the operator $W(c)$ is
invertible from the right. If $s\in\ker W(d)$, then the element
$(W_r^{-1}(c) W(\widetilde{a}^{-1})s,s)^T\in\ker W(V(a,b)))$. By
Lemma \ref{l1}, the element
\[
2\vp^{\pm}(s):= W_r^{-1}(c)W(\widetilde{a}^{-1})s\mp JQW^0(c)P
W_r^{-1}(c)W(\widetilde{a}^{-1})s\pm JQW^0(\widetilde{a}^{-1}) s
\]
belongs to the null space $\ker(W(a)\pm H(b))$ of the corresponding
operator $W(a)\pm H(b)$.

\begin{lem}\label{l2}
Let $(c,d)$ be the subordinated pair for a matching pair $(a,b)\in
G\times G$. If the operator $W(c)$ is right-invertible, then for
every $s\in\ker W(d)$ the following relations
\[
(W(\widetilde{b})+H(\widetilde{a}))\vp_+(s)=\mathbf{P}^+(d)s,\quad
(W(\widetilde{b})-H(\widetilde{a}))\vp_-(s)=\mathbf{P}^-(d)s,
\]
hold. Thus the corresponding mappings $\vp_+: \im \mathbf{P}^+(d)
\to \im \mathbf{P}^+(d) $ and $\vp_-: \im\mathbf{P}^-(d)\to
\im\mathbf{P}^-(d)$, are injective operators.
\end{lem}

\textbf{Proof.} Assuming that $s\in\ker W(d)$, one can show that the
operator $W(\widetilde{b})+H(\widetilde{a})$ sends $\vp^+(s)$ into
$\mathbf{P}^+(d)s$ and the operator
$W(\widetilde{b})-H(\widetilde{a})$ sends $\vp^-(s)$ into
$\mathbf{P}^-(d)s$. The proof of these facts is based on relations
\eqref{cst3} and runs similarly to the proof of \cite[Lemma
3.6]{DS:2013e}. \rbx

\begin{prop}\label{p2}
Let $(c,d)$ be the subordinated pair for a matching pair $(a,b)\in G
\times G$. If the operator $W(c)$ is right-invertible, then
\begin{align*}
\ker(W(a)+H(b))&=\vp^+(\im \mathbf{P}^+(d)) \dotplus\im \mathbf{P}^-(c), \\
\ker(W(a)-H(b))&=\vp^-(\im \mathbf{P}^-(d)) \dotplus\im
\mathbf{P}^+(c).
\end{align*}
\end{prop}

\textbf{Proof.} Using the invertibility of the operator $E_1$ and
Proposition \ref{p3.4}, one obtains
\[
\ker\diag (W(a)+H(b),W(a)-H(b))=E_1(\widehat{\Omega}(d))\dotplus
E_1(\Omega(c)).
\]
Apparently,
$\widehat{\Omega}(d)=\widehat{\Omega}_+(d)\dotplus\widehat{\Omega}_-(d)$,
$\Omega(c)=\Omega_+(c)\dotplus\Omega_-(c)$, where
\begin{align*}
\widehat{\Omega}_+(d)&=\big\{(W_r^{-1}(c)W(\widetilde{a}^{-1})s,s)^T:s\in\im\mathbf{P}^+(d)\big\},\\
\widehat{\Omega}_-(d)&=\big\{(W_r^{-1}(c)W(\widetilde{a}^{-1})s,s)^T:s\in\im\mathbf{P}^-(d)\big\},\\
\Omega_+(c)&=\big\{(s,0)^T:s\in\im\mathbf{P}^+(c)\big\},\\
\Omega_-(c)&=\big\{(s,0)^T:s\in\im\mathbf{P}^-(c)\big\}.
\end{align*}
Hence
\begin{multline*}
\ker\diag (W(a)+H(b),W(a)-H(b))\\
=E_1(\widehat{\Omega}_+(d))\dotplus
E_1(\widehat{\Omega}_-(d))\dotplus E_1(\Omega_+(c))\dotplus
E_1(\Omega_-(c)) =E_1(\ker W(V(a,b)).
\end{multline*}
It is clear that if $\phi\in\ker(W(a)+H(b))$, then
$(\phi,0)^T\in\ker\diag (W(a)+H(b),W(a)-H(b)).$ Now we want to find
that uniquely defined element $(\alpha,\beta)^T$ from the kernel of
the operator $W(V(a,b))$, which is sent into the element $(\phi,0)$
by the operator $E_1$. It can be uniquely represented in the form
\[
(\alpha,\beta)^T\!\!=
 (W_r^{-1}(c)W(\widetilde{a}^{-1})s_+,s_+)^T\!\!
+(W_r^{-1}(c)W(\widetilde{a}^{-1})s_-,s_-)^T\!\!
+(v_+,0)^T\!\!+(v_-,0)^T,
\]
where $s_{\pm}\in\im\mathbf{P}^{\pm}(d)$,
$v_{\pm}\in\im\mathbf{P}^{\pm}(c)$. Then
\begin{align*}
(\phi,0)^T&=E_1((\alpha,\beta)^T) \\
&=E_1((W_r^{-1}(c)W(\widetilde{a}^{-1})s_+,s_+)^T)
 +E_1((W_r^{-1}(c)W(\widetilde{a}^{-1})s_-,s_-)^T) \\
&\quad+E_1((v_+,0)^T)+E_1((v_-,0)^T) \\
&=(\vp_+(s_+),\vp_-(s_+))^T+(\vp_+(s_-),\vp_-(s_-))^T+(0,v_+)^T+(v_-,0)^T.
\end{align*}
Thus
\[
\phi=\vp_+(s_+)+\vp_+(s_-)+v_-, \quad 0=\vp_-(s_+)+\vp_-(s_-)+v_+.
\]
However, since $\vp_+(s_-)\in\ker (W(a)+H(b))$ and
$E_2((\vp_+(s_-),0)^T)=(\vp_+(s_-),0)^T$ according to Lemma
\ref{l2}, we get $\vp_+(s_-)\in \im\mathbf{P}^-(c)$. Analogously,
one can show that $\vp_-(s_+)\in \mathbf{P}^+(c)$. It implies that
$\vp_-(s_-)=-(\vp_-(s_+)+v_+)\in\im\mathbf{P}^+(c) $ and
$E_2((0,\vp_-(s_-))^T)=(\vp_-(s_-),0)^T$ because
$\vp_-(s_-)\in\im\mathbf{P}^+(c) $. On the other hand, due to Lemma
\ref{l2}, one has $E_2((0,\vp_-(s_-))^T)=(\vp_-(s_-),s_-)^T$. The
comparison of the two expressions for the element
$E_2((0,\vp_-(s_-))^T)$ gives $s_-=0$, and therefore $\vp_-(s_-)=0$.
This implies $\vp_-(s_+)=-v_+$. Consequently,
\[
(\alpha,\beta)^T=(W_r^{-1}(c)W(\widetilde{a}^{-1})s_+,s_+)^T-(\vp_-(s_+),0)^T+(v_-,0)^T,
\]
which leads to the relation
\[
E_1((\alpha,\beta)^T)=(\vp_+(s_+),\vp_-(s_+))^T-(0,\vp_-(s_+))^T+(v_-,0)^T=(\vp_+(s_+)+v_-,0)^T.
\]
Thus $\vp_+(s_+)+v_-\in \ker(W(a)+H(b))$. This result shows that
$\ker(W(a)+H(b))$ is the sum of its subspaces
$\vp_+(\im\mathbf{P}^+(d))$ and $\im \mathbf{P}^-(c)$. Recalling
that
$(W_r^{-1}(c)W(\widetilde{a}^{-1})s_+,s_+)^T-(\vp_-(s_+),0)^T\in\widehat{\Omega}_+(d)\dotplus\Omega_+(c)$
and $(v_-,0)^T \in\Omega_-(c)$, one finally obtains
\[
\ker(W(a)+H(b))=\vp_+(\im
\mathbf{P}^+(d))\dotplus\im\mathbf{P}^-(c).
\]
The relation
\[
\ker(W(a)-H(b))=\vp_-(\im \mathbf{P}^-(d))\dotplus\im\mathbf{P}^+(c)
\]
can be verified analogously. \rbx

\begin{cor}\label{cc1}
Let $(c,d)$ be the subordinated pair for a matching pair $(a,b)\in
G\times G$ satisfying the conditions of Proposition \ref{c2}. Then
\begin{align*}
\dim\ker(W(a)+H(b))&=\dim\im\mathbf{P}^+(d)+\dim\im\mathbf{P}^-(c),\\
\dim\ker(W(a)-H(b))&=\dim\im\mathbf{P}^-(d)+\dim\im\mathbf{P}^+(c).
\end{align*}
\end{cor}

\begin{rem}
If $(a,b)\in G\times G$ is a Fredholm matching pair, i.e., if
$W(c),W(d)$ are Fredholm operators, then $W(a)\pm H(b)$ are Fredholm
operators and
\begin{equation}\label{cst15}
\ind(W(a)+H(b))+\ind (W(a)-H(b))=\ind W(c)+\ind W(d).
\end{equation}
We conjecture that if one of the operators $W(a)+H(b)$ or
$W(a)-H(b)$ is Fredholm, then so is the other and relation
\eqref{cst15} holds.
\end{rem}

\section{Kernels of Wiener--Hopf plus Hankel operators. Specification}\label{s3}

In this section we study the kernels of Wiener--Hopf plus Hankel
operators $W(a)+H(b)$ in the case where the generating functions
$a,b\in G$ satisfy matching condition \eqref{cst9} and $W(c),W(d)$
are mainly Fredholm operators such that
\[
0\leq |\ind W(c)|,|\ind W(d)|\leq 1.
\]
Recall that $a$ is supposed to be invertible in $G$. In view of
Theorem \ref{t1}, on has $\nu(c)=\nu(d)=0$ and
$0\leq|n(c)|,|n(d)|\leq 1$.

In order to formulate our first result we need the following lemma.

\begin{lem}\label{l3}
 Let $\chi(t):=(t-i)/(t+i)$, $t\in\sR$.
\begin{enumerate}
\item If the function $\psi$ is defined by
\[
\psi(t):=\left\{\begin{array}{ll}e^{-t} &\text{ if }t>0, \\0&\text{
if } t<0,\end{array} \right.
\]
then $W^0(\chi^{-1})\psi=-\widetilde{\psi}$.

\item On each space $L^p(\sR^+)$, $1\leq p\leq\infty$
the operator $W(\chi^{-1})$ has a one-dimensional kernel generated
by the function $\psi_0(t)=e^{-t}$, $t>0$.
\end{enumerate}
\end{lem}

\textbf{Proof.} Assertion (i) can be obtained by using the relation
\[
(W^0(\chi^{-1})g)(t)=g(t)-2\int_t^\infty e^{t-s} g(s)\,ds, \quad
-\infty<t<\infty,
\]
which is valid for all $g\in L^p(\sR)$, $1\leq p\leq\infty$,
\cite{GF1974}.

Assertion (ii) is well known. It can be proved by the
differentiation of the identity
\[
\vp(t)=2 \int_t^\infty e^{t-s}\vp(s)\,ds.
\]
Moreover, one has
\[
(W(\chi^{-1})g)(t)=g(t)-2\int_t^\infty e^{t-s} g(s)\,ds, \quad
0<t<\infty,
\]
(see \cite{GF1974}). Note that assertion~(ii) also follows from
assertion~(i). \rbx

Now we can derive the following version of the Coburn--Simonenko
Theorem.

\begin{thm}\label{t2}
Let $a\in G$ be invertible and let $A$ denote any of the four
operators $W(a)-H(a\chi)$, $W(a)+H(a\chi^{-1})$, $W(a)\pm H(a)$.
Then at least one of the spaces $\ker A$ or $\coker A$ is trivial.
\end{thm}

\textbf{Proof.}  \textsf{Part 1:} Let us start with the operator
$W(a)+H(a\chi)$. The function $\chi$ satisfies the relation
$\widetilde{\chi}=\chi^{-1}$, so the duo $(a,a\chi)$ is a matching
pair with the subordinated pair $(c,d)$ with $c=\chi^{-1}$ and $d=a
\widetilde{a}^{-1} \chi$. Moreover, the operator $W(\chi^{-1})$ is
invertible from the right and one of its right inverses is the
operator $W(\chi)$. Thus the theory of Section \ref{s2} applies. As
it was pointed out earlier, the kernel of this operator is $ \ker
W(\chi^{-1})=\{\mathbf{c}\psi_0:\mathbf{c}\in\sC\}, $ where
$\psi_0(t)=e^{-t}$, $t>0$. In order to apply Proposition \ref{p2} we
have to identify, in particular, the projections
$\mathbf{P}^{\pm}(\chi^{-1})$ acting on the space $\ker
W(\chi^{-1})$. But $\mathbf{P}^+(\chi^{-1})$ and
$\mathbf{P}^-(\chi^{-1})$ are complimentary projections on the
one-dimensional space $\ker W(\chi^{-1})$. Therefore, one of these
projections is just the identity operator whereas the other one is
the zero operator. Consider next the expression
$JQW^0(\chi^{-1})P\psi_0$. By Lemma \ref{l3}(i) one has
\[
JQW^0(\chi^{-1}) P\psi_0=JQW^0(\chi^{-1})\psi=-JQ\widetilde{\psi},
\]
so that $JQW^0(\chi^{-1}) P\psi_0=-\psi_0$ and
$\mathbf{P}^-(\chi^{-1})=I$ on $\ker W(\chi^{-1})$.

According to Proposition \ref{p2}, the kernels of the operators
$W(a)+H(a\chi)$ and $W(a)-H(a\chi)$ can be represented in the form
\begin{equation}\label{cst16}
\begin{aligned}
\ker(W(a)-H(a\chi))&=\vp^-(\im \mathbf{P}^-(d)) , \\
\ker(W(a)+H(a\chi))&=\vp^+(\im \mathbf{P}^+(d))
\dotplus\{\mathbf{c}\psi_0: \mathbf{c}\in \sC\}.
\end{aligned}
\end{equation}
If $\dim\ker W(d)>0$, then $\coker (W(a)\pm H(a\chi))=\{0\}$.
Indeed, relation \eqref{cst11} and the familiar Coburn--Simonenko
Theorem for the operator $W(d)$ show that $\coker
W(V(a,a\chi))=\{0\}$. Taking into account representation
\eqref{cst8}, one obtains that the cokernel of each of the operators
$W(a)+H(a\chi)$ and $W(a)-H(a\chi)$ contains the zero element only.

Let us now assume that $\ker W(d)=\{0\}$. Then the first relation
\eqref{cst16} implies that $\ker(W(a)-H(a\chi))=0$. Hence, the
operator $W(a)-H(a\chi)$ is subject to Coburn--Simonenko Theorem.

\smallskip
\textsf{Part 2:} Consider the operator $W(a)+H(a\chi^{-1})$ and note
that $W(c)=W(\chi)$ is not right-invertible, so that Proposition
\ref{p2} cannot be directly used in this situation. Nevertheless,
the case at hand can be reduced to the operators studied. Thus the
operators $W(a)\pm H(a\chi^{-1})$ can be represented in the form
\begin{equation}\label{cst17}
W(a)\pm H(a\chi^{-1})=(W(a\chi^{-1})\pm H(a\chi^{-1}\chi))W(\chi).
\end{equation}
The proof of \eqref{cst17} follows from \eqref{cst4} and relation
$H(\chi)W(\chi)=0$. Setting $\alpha:=a\chi^{-1}$, we get
\begin{equation}\label{cst18}
W(a)\pm H(a\chi^{-1})=(W(\alpha)\pm H(\alpha\chi)) W(\chi).
\end{equation}
The operators of the form $W(\alpha)\pm H(\alpha\chi)$ in the
right-hand side of \eqref{cst17} have been just studied, and we
already know that the function $\psi_0$ belongs to the kernels of
both operators $W(\alpha)+H(\alpha\chi)$ and $W(\chi^{-1})$. Since
$W(\chi^{-1})W(\chi)=I$ it follows that $\psi_0\notin \im W(\chi)$.
Consider now the projection $Q_0:=W(\chi)W(\chi^{-1})$ which
projects the space $L^p(\sR^+)$, $1\leq p\leq \infty$ onto $\im
W(\chi)$ parallel to $\ker W(\chi^{-1})$.

Assume first that $\ker W(d)=\{0\}$ and note that for the matching
pairs $(a,a\chi^{-1})$ and $(\alpha,\alpha\chi)$, the corresponding
subordinated pairs $(c,d)$ have the same element $d$, namely, $d=a
\widetilde{a}^{-1}\chi$. Then \eqref{cst18} shows that
$\ker(W(a)+H(a\chi^{-1}))=\{0\}$. Further, if $\dim\ker W(d)>0$,
then the space $\ker(W(\alpha)+H(\alpha \chi))$ decomposes as
follows
\[
\ker(W(\alpha)+H(\alpha \chi))=\ker W(\chi^{-1})\oplus
Q_0(\ker(W(\alpha)+H(\alpha \chi))).
\]
However, as was already shown, the operator
$W(\alpha)-H(\alpha\chi)$ is right-invertible and
\[
\ker W(\chi^{-1})\subset\ker W(\alpha)+H(\alpha\chi).
\]
Therefore, relation \eqref{cst17} implies that the operator
$W(a)+H(a\chi^{-1})$ maps $L^p(\sR^+)$ onto $L^p(\sR^+)$, so it is
subject to the Coburn--Simonenko Theorem.

\smallskip
\textsf{Part 3:} It remains to consider the operators $W(a)\pm
H(a)$. For these operators the element $c$ in the corresponding
subordinated pair is either $1$ or $-1$, and our claim follows
immediately from the Coburn--Simonenko Theorem for scalar
Wiener--Hopf operators and from relations \eqref{cst11} and
\eqref{cst8}. \rbx

\begin{rem}\label{r5}
The proof of Theorem \ref{t2} shows that this theorem remains true
for more general generating functions, for instance, in the case
where $a$ and $b$ belong to the algebras $G_p$, $1\leq p\leq\infty$
studied in \cite[Chapter VII]{GF1974}.
\end{rem}

The reader can also observe that, in fact, we have proved a bit more
than Theorem \ref{t2} states. A more detailed result can be
formulated as follows.

\begin{cor}\label{c2}
Let $a\in G$  be invertible. Then
\begin{enumerate}
\item If $\dim\ker W(d)=0$, then
\begin{align*}
\ker(W(a)-H(a\chi))=\{0\},\quad
\ker(W(a)+H(a\chi))=\{\mathbf{c}\psi_0: \mathbf{c}\in\sC\},
\end{align*}
and if $\dim\ker W(d)>0$, then $\coker(W(a)\pm H(a\chi))=\{0\}$.
\item If $\dim\ker W(d)=0$, then $\ker(W(a)\pm H(a\chi^{-1}))=\{0\}$,\\
and if $\dim\ker W(d)>0$, then $\coker(W(a)+H(a\chi^{-1}))=\{0\}$.
\end{enumerate}
\end{cor}

Let us emphasize that the description of the projections
$\mathbf{P}^{\pm}(\chi^{-1})$ did play an important role in our
considerations. In the general case one has to study the projections
$\mathbf{P}^{\pm}(g)$ for the functions $g$ satisfying the relation
$g \widetilde{g}=1$. Because of the space restriction, we are not
going to pursue this matter here. Nevertheless, let us consider the
case where $\nu(g)=0$ and $n(g)=-1$, which is one of the simplest
generalization of the situation $g=\chi^{-1}$. In order to handle
this case we need a result from \cite[Chapter VII]{GF1974}.

\begin{prop}\label{p5}
Each invertible function $g\in G$ admits the factorization of the
form
\begin{equation}\label{cst19}
g(t)=g_-(t)e^{i\nu t}\Big(\frac{t-i}{t+i}\Big)^n g_+(t), \quad
-\infty<t<\infty,
\end{equation}
where $g_+^{\pm1}\in G^+$, $g_-^{\pm1}\in G^-$, $\nu=\nu(g)$ and
$n=n(g)$. Moreover, under the agreement $g_-(0)=1$, the
factorization factors $g_+$ and $g_-$ are uniquely defined.
\end{prop}

Note that the proof of Theorem \ref{t1} is based on Proposition
\ref{p5}.

\begin{defn}\label{def1}
Suppose that $g\in G$ satisfies the condition $g \widetilde{g}=1$
and set
\[
\xi(g)=(-1)^n g(0), \quad n=n(g).
\]
\end{defn}

\begin{thm}\label{t3}
If $g\in G$ and $g \widetilde{g}=1$, then $\xi(g)=\pm1$ and the
factorization~\eqref{cst19} takes the form
\begin{equation}\label{cst20}
g(t)=\left(\xi(g)\,\widetilde{g}_+^{-1}(t)\right)e^{i\nu
t}\Big(\frac{t-i}{t+i}\Big)^n g_+(t)
\end{equation}
with $\widetilde{g}_+^{\pm1}(t)\in G^-$ and
$g_-(t)=\xi(g)\,\widetilde{g}_+^{-1}(t)$.
\end{thm}

\textbf{Proof.} Using the condition $g^{-1}=\widetilde{g}$, we get
from \eqref{cst19} that
\[
g_+^{-1}(t) e^{-i\nu t}\Big(\frac{t-i}{t+i}\Big)^{-n}
g_-^{-1}(t)=\widetilde{g}_-(t) e^{-i\nu
t}\Big(\frac{t-i}{t+i}\Big)^{-n}\widetilde{g}_+(t),
\]
where $\nu=\nu(g)$, $n=n(g)$.

Note that $\widetilde{g}^{\pm1}_-\in G^+$,
$\widetilde{g}^{\pm1}_+\in G^-$, as easy computations show.
Therefore,
\[
g_+^{-1}\widetilde{g}_-^{-1}=g_-\widetilde{g}_+,
\]
and $g_+^{-1}\widetilde{g}_-^{-1}\in G^+$, $g_-\widetilde{g}_+\in
G^-$. It follows that there is a constant $\xi\in\sC$ such that
$g_+^{-1}\widetilde{g}_-^{-1}=\xi=g_-\widetilde{g}_+$, and
$g_-=\xi\, \widetilde{g}_+^{-1}$. For the function $g_0=g_+g_-$ we
have $g_0 \widetilde{g}_0=1$. Therefore,
\[
1=g_0 \widetilde{g}_0=(\xi
g_+\widetilde{g}_+^{-1})(\xi\widetilde{g}_+g_+^{-1})=\xi^2.
\]
For $t=0$, which is one of the fixed points of the operator $J$, the
equation $g_0=\xi g_+\widetilde{g}_+^{-1}$ implies $g_0(0)=\xi$, and
$g_0(0)=g(0)(-1)^n$ (see \eqref{cst20}). Thus we obtain that
$\xi=g(0)(-1)^n$ which completes the proof. \rbx

Now we again use the notation
\[
\chi^{\pm1}(t)=\Big(\frac{t-i}{t+i}\Big)^{\pm1}, \quad t\in\sR.
\]

\begin{thm}\label{t4}
Let $g\in G$, $g \widetilde{g}=1$, $\nu(g)=0$ and $n(g)=-1$. Then
\[
\im
\mathbf{P}^{\pm}(g)=\bigg\{\mathbf{c}\Big(\frac{1\mp\xi(g)}{2}\Big)W(g_+^{-1})\psi_0:\mathbf{c}\in\sC\bigg\}.
\]
\end{thm}

\textbf{Proof.} It is easily seen that $\ker
W(g)=\left\{\mathbf{c}W(g_+^{-1})\psi_0:\mathbf{c}\in\sC\right\}$,
According to the definition of projections $\mathbf{P}^{\pm}(g)$ we
have to compute the expression
\[
JQW^0(g)PW(g_+^{-1})\psi_0.
\]
We have
\begin{align*}
JQW^0(g)PW(g_+^{-1})&=JQW^0(g_-)W^0(\chi^{-1})W^0(g_+)W^0(g_+^{-1})P \\
&=JQW^0(g_-)W^0(\chi^{-1})P.
\end{align*}
Recall that by Lemma \ref{l3},
$W^0(\chi^{-1})P\psi_0=W^0(\chi^{-1})\psi=-\widetilde{\psi}$, and
using Theorem~\ref{t3} we get
\begin{align*}
JQW^0(g)PW(g_+^{-1})\psi_0&=-JQ W^0(g_-)\widetilde{\psi}=-W^0(\widetilde{g}_-)\psi\\
&=-P\xi(g)W^0(g_+^{-1})P\psi_0=-\xi(g)W^0(g_+^{-1})\psi_0,
\end{align*}
and we are done. \rbx

The next result is a generalization of Theorem \ref{t2}.

\begin{thm}\label{t5}
Let $a,b\in G$ constitute a matching pair, $a$ be invertible in $G$
and let $(c,d)$ be the subordinated pair for $(a,b)$. If $A$ denotes
one of the following operators
\begin{enumerate}
\item $W(a)\pm H(b)$ with $\nu(c)=0$, $n(c)= 1$ and $\xi(c)=\pm1$;
\item $W(a)\mp H(b)$ with $\nu(c)=0$, $n(c)=-1$ and $\xi(c)=\pm1$;
\item $W(a)\pm H(b)$ with $\nu(c)=0$ and $n(c)=0$
\end{enumerate}
considered on the space $L^p(\sR^+)$, then at least one of the
spaces $\ker A$ or $\coker A$ is trivial.
\end{thm}

\textbf{Proof.} The proof mimics that of Theorem \ref{t2} with minor
modifications. First, we note that the case $\xi(c)=-1$ can be
reduced to the case $\xi(c)=1$ via rearrangements
$W(a)+H(b)=W(a)-H(-b)$ and $W(a)-H(b)=W(a)+H(-b)$. Therefore, we
only consider the situation $\xi(c)=1$ in the cases (i) and (ii).
Further, one has to use Theorem \ref{t4} instead of the description
of the projections $\mathbf{P}^{\pm}(\chi^{\pm1})$. Consider the
operator $W(a)+H(b)$ in the case where $\nu(c)=0$ and $n(c)=1$.
Representing the operator $W(a)\pm H(b)$ in the form
\[
W(a)\pm H(b)=(W(a\chi^{-1})\pm H(b\chi))W(\chi),
\]
we observe that $(a\chi^{-1},b\chi)$ is a matching pair with the
subordinated pair $(c\chi^{-2}, d)$ and $\ind W(c\chi^{-2})=-1$,
$\im\mathbf{P}^+(c\chi^{-2})=\ker
W(c\chi^{-2})=\{\mathbf{c}W(c_+^{-1})\psi_0:c\in\sC\}$. Let us also
note that $\ker W(c\chi^{-2})=\ker W(c_+\chi^{-1})$ and
$W(c_+\chi^{-1})\,W(c_+^{-1}\chi)=I$. Hence, $\ker W(c\chi^{-2})
\cap \im W(c_+^{-1}\chi)=\{0\}$. Since obviously $\im
W(c_+^{-1}\chi)=\im W(\chi)$, we obtain $\ker W(c\chi^{-1})\cap \im
W(\chi)=\{0\}$.

Now one can proceed similarly to Part 2 in the proof of Theorem
\ref{t2}. \rbx

\begin{cor}\label{c3}
Assume that $a,b\in G$ constitute a matching pair with the
subordinated pair $(c,d)$ such that $\xi(c)=1$. Then
\begin{enumerate}
\item If $\dim\ker W(d)=0$, and $\ind W(c)=1$, then
\begin{align*}
\ker(W(a)-H(b))=\{0\},\quad
\ker(W(a)+H(b))=\{\mathbf{c}W(c_+^{-1})\psi_0: \mathbf{c}\in\sC\},
\end{align*}
and if $\dim\ker W(d)>0$, then $\coker(W(a)\pm H(b))=\{0\}$.
\item If $\dim\ker W(d)=0$, and $\ind W(c)=-1$, then $\ker(W(a)\pm H(b))=\{0\}$,\\
and if $\dim\ker W(d)>0$, then $\coker(W(a)+H(b))=\{0\}$.
\end{enumerate}
\end{cor}

An interesting and important subclass of the operators considered in
this paper comprises the identity plus Hankel operators. Let us
specify the above results in this situation

\begin{cor}\label{c4}
If $b\in G$ is a matching function, then $(1,b)$ is a matching pair
with the subordinated pair $(\widetilde{b},b)$, and if $A$ denotes
any of the operators
\begin{enumerate}
\item $I-H(b)$ with $\nu(\widetilde{b})=0$, $n(\widetilde{b})=-1$ and $\xi(\widetilde{b})=1$;
\item $I+H(b)$ with $\nu(\widetilde{b})=0$, $n(\widetilde{b})= 1$ and $\xi(\widetilde{b})=1$;
\item $I\pm H(b)$ with $\nu(\widetilde{b})=0$ and $n(\widetilde{b})=0$,
\end{enumerate}
considered on the space $L^p(\sR^+)$, then $\ker A$ or $\coker A$ is
trivial.
\end{cor}

Now we revisit Theorem \ref{t2} and consider the operators $W(a)\pm
H(a\chi)$ and $W(a)\pm H(a\chi^{-1})$ under additional assumptions.

\smallskip
$\mathbf{1^0}$. Suppose that $\nu(a)=n(a)=0$ and $b=a\chi$. The
subordinated pair $(c,d)$ is given by the elements $c=\chi^{-1}$ and
$d=a \widetilde{a}^{-1} \chi$. Thus
\[
\ind W(c)=1, \quad \ind W(d)=-1, \quad \xi(c)=\xi(d)=1.
\]
According to \eqref{cst15} we have
\begin{equation}\label{cst21}
\ind(W(a)+H(a\chi))+\ind (W(a)-H(a\chi))=0.
\end{equation}
Further, by Corollary \ref{c2}(i) we also have
\[
\ker(W(a)-H(a\chi))=0, \quad \ker
(W(a)+H(a\chi))=\{\mathbf{c}\psi_0: \mathbf{c}\in \sC \}.
\]
In order to describe the cokernels of the above operators we make
use of the adjoint operators. If $p\in [1,\infty)$, then according
to \eqref{cst10} the adjoint operators have the form
$W(\overline{a})\pm H(\overline{\widetilde{a}}\chi)$, and the duo
$(\overline{a},\overline{\widetilde{a}}\chi)$ is a matching pair
with the subordinated pair $(\overline{d},\overline{c})$, so that
$\ind W(\overline{d})=1, \ind W(\overline{c})=-1$ and
$\xi(\overline{d})=1$. By Corollary \ref{c3}(ii),
$\ker(W(\overline{a})-H(\overline{\widetilde{a}}\chi))=\{0\}$, which
finally proves that the operator $W(a)-H(a\chi)$ is invertible. Note
that this result is also true for the space $L^\infty(\sR^+)$.
Indeed, the operator
$W(\overline{a})-H(\overline{\widetilde{a}}\chi)$ acts on the space
$L^1(\sR^+)$ and the above considerations show that
$\dim\ker(W(\overline{a})-H(\overline{\widetilde{a}}\chi))=0$. The
adjoint of this operator acts on the space $L^\infty(\sR^+)$ and is
equal to the operator $W(a)-H(a\chi)$, the kernel of which is
trivial. Therefore, the operator
$W(\overline{a})-H(\overline{\widetilde{a}}\chi)$ is invertible on
the space $L^1(\sR^+)$. Consequently, its adjoint $W(a)-H(a\chi)$ is
invertible on $L^\infty(\sR^+)$. Then relation \eqref{cst21}
immediately implies that $\ind (W(a)+H(a\chi))=0$. Note that the
operator $W(a)+H(a\chi)$ provides an example of operators where both
spaces $\ker(W(a)+H(a\chi))$ and $\coker(W(a)+H(a\chi))$ are
nontrivial.

\smallskip
$\mathbf{2^0}$. Suppose that $\nu(a)=0$, $n(a)=-1$ and $b=a\chi$.
For the subordinated pair $(c,d)$ we have $c=\chi^{-1}$ and
$d=a\widetilde{a}^{-1} \chi$ so that $\ind W(c)=1$, $\ind W(d)=1$,
$\xi(d)=1$. Since $\ind W(d)=1$, Corollary \ref{c2}(i) indicates
that $\coker(W(a)\pm H(a\chi))=\{0\}$. besides, $\dim\ker (W(a)\pm
H(a\chi))=1$ by Proposition \ref{p2}.

\smallskip
$\mathbf{3^0}$. Suppose that $\nu(a)=n(a)=0$ and $b=a\chi^{-1}$.
Since $c=\chi$, the operator $W(c)$ is not invertible from the
right. Write
\begin{equation}\label{cst22}
W(a)\pm H(a\chi^{-1})=(W(a\chi^{-1})\pm H(a\chi^{-1}\chi))W(\chi),
\end{equation}
and set $\alpha:=a\chi^{-1}$. The operators $W(\alpha)\pm
H(\alpha\chi)$ are considered in $\mathbf{2^0}$, so we have
\begin{align*}
\dim\ker(W(\alpha)\pm H(\alpha\chi))&=1\\
\dim\coker (W(\alpha)\pm H(\alpha\chi))&=0.
\end{align*}
According to the Part 2 in the proof of Theorem \ref{t2}, one has
$\ker(W(a)+H(a\chi^{-1}))=\{0\}$. This and the relation
$\dim\coker(W(\alpha)+H(\alpha\chi))=0$ show the invertibility of
the operator $W(a)+H(a\chi^{-1})$. Due to Proposition \ref{p2} (see
also \eqref{cst16}) we know that the kernel of the operator
$W(\alpha)-H(\alpha\chi)$ is spanned on the element
\begin{multline}\label{cst23}
\kappa= W(\chi) W(\widetilde{\alpha}^{-1})W(d_+^{-1})\psi_0
+JQW^0(\chi^{-1})PW^0(\chi)PW(\widetilde{\alpha}^{-1})W(d_+^{-1})\psi_0 \\
-JQW^0(\widetilde{\alpha}^{-1}) P W(d_+^{-1}) \psi_0,
\end{multline}
where we used the fact that $W(\chi)$ is a right inverse for the
operator $W(\chi^{-1})$ and where $d_+^{-1}$ arises from the
factorization \eqref{cst20} of the function
$d=a\widetilde{a}^{-1}\chi^{-1}$. Note that the first term in
\eqref{cst23} belongs to the set $\im W(\chi)$, whereas the second
one is equal to zero. Thus the operator $W(a)-H(a\chi^{-1})$ is
invertible if and only if $H(\alpha^{-1}) W(d_+^{-1})\psi_0\notin\im
W(\chi)$. On the other hand, if this condition is not satisfied, the
operator $W(a)-H(a\chi^{-1})$ presents an example of a Wiener--Hopf
plus Hankel operator with one-dimensional kernel and cokernel.

\smallskip
$\mathbf{4^0}$. Suppose that $\nu(a)=0$, $n(a)=1$ and
$b=a\chi^{-1}$. Let us use representation \eqref{cst22} and set
$\alpha=a\chi^{-1}$. It follows from Part $\mathbf{1^0}$ that
$W(\alpha)-H(\alpha\chi)$ is invertible whereas the operator
$W(\alpha)+H(\alpha\chi)$ has one-dimensional kernel and cokernel.
Since
\[
\ker
(W(\alpha)+H(\alpha\chi))=\{\mathbf{c}\psi_0:\mathbf{c}\in\sC\}\cap\im
W(\chi)=\{0\},
\]
we conclude that the operator $W(a)+H(a\chi^{-1})$ has trivial
kernel and a cokernel of dimension $1$. Of course, the same
conclusion is valid for the operator $W(a)-H(a\chi)$.

It is worth noting that a similar consideration with natural
amendments can be used in the contest of Theorem \ref{t5}. Let us
restrict ourselves to the operators $I+H(b)$ with the generating
function $b$ satisfying the condition $b \widetilde{b}=1$. Then
$(1,b)$ is a matching pair with the subordinated pair
$(\widetilde{b},b)$.

\smallskip
$\mathbf{5^0}$. Suppose that $\nu(b)=n(b)=0$. Then the operators
$W(b)$ and $W(\widetilde{b})$ are invertible and relations
\eqref{cst8}, \eqref{cst11} already show that $I+H(b)$ and $I-H(b)$
are invertible operators.

Assume next that $\nu(b)=0$ but $n(b)=1$ and $\xi(\widetilde{b})=1$.
Then $\ind W(\widetilde{b})=1$ and $\ind W(b)=-1$. By Corollary
\ref{c3}(i), one has
\[
\ker (I-H(b))=\{0\}, \quad
\ker(I+H(b))=\{\mathbf{c}W(b_+)\psi_0:\mathbf{c}\in\sC\}.
\]
Similarly to Part $\mathbf{1^0}$ one shows that the operator
$I-H(b)$ is invertible and $\ind(I+H(b))=0$.

Finally, let us assume that $\nu(b)=0$, $n(b)=-1$ and
$\xi(\widetilde{b})=1$. Since $\ind W(\widetilde{b})=-1$, we will
use the relation $I\pm H(b)=(W(\chi^{-1})\pm H(b\chi)) W(\chi)$. It
is clear that $(\chi^{-1},b\chi)$ is a matching pair with the
subordinated pair $(\widetilde{b}\chi^{-2}, b)$ and $\ind
W(\widetilde{b}\chi^{-2})\!=\ind W(b)=1$. Analogously to Part
$\mathbf{2^0}$ we obtain that
\[
\coker (W(\chi^{-1})\pm H(b\chi))=\{0\}.
\]
Moreover, by Proposition \ref{p2}, $\dim\ker (W(\chi^{-1})\pm
H(b\chi))=1$ and since
\[
\ker
(W(\chi^{-1})+H(b\chi))=\{\mathbf{c}W(b_+)\psi_0:\mathbf{c}\in\sC\}\cap\im
W(\chi)=\{0\},
\]
the operator $I+H(b)$ is invertible. If $\ker
W(\chi^{-1})-H(b\chi)\cap \im W(\chi)=\{0\}$, then $I-H(b)$ is
invertible. Otherwise, $\ind(I-H(b))=0$, but this operator is not
invertible.

%%%%%%%
%\bibliographystyle{amsplain}
 % \bibliographystyle{acm}
% \bibliography{E:/TeXX/JabRef/DWP_2010}

\begin{thebibliography}{99} %% n is number of items, or the largest label


\bibitem{BE2004}
E.\,L. Basor, T. Ehrhardt, \emph{Factorization theory for a class
  of {T}oeplitz {$+$} {H}ankel operators}, J. Operator Theory
  \textbf{51}(2)(2004), 411--433.

\bibitem{BE2006}
E.\,L. Basor, T. Ehrhardt, \emph{Factorization of a class of
{T}oeplitz + {H}ankel operators and the {$A_p$}-condition}, J.
Operator Theory \textbf{55}(2) (2006), 269--283.

\bibitem{BE:2013}
E.\,L. Basor, T. Ehrhardt, \emph{Fredholm and invertibility theory
for a special class of {T}oeplitz + {H}ankel operators}, J. Spectral
Theory \textbf{3}(3) (2013), 171--214.

\bibitem{BS:2006}
A. B{\"o}ttcher, B. Silbermann, \emph{Analysis of {T}oeplitz
operators}, second ed., Springer Monographs in Mathematics,
Springer-Verlag, Berlin, 2006.

\bibitem{CS:2012}
L.\,P. Castro, A.\, S. Silva, \emph{Defect numbers of singular
integral operators with {C}arleman shift and almost periodic
coefficients}, J. Math. Anal. Appl. \textbf{387}(1) (2012), 66--76.

\bibitem{DS:2013d}
V.\,D. Didenko, B.~Silbermann, \emph{Some results on the
invertibility of {T}oeplitz plus {H}ankel operators}, Ann. Acad.
Sci. Fenn. Math. \textbf{39} (2014), 443--461.


\bibitem{DS:2013e}
V.\,D. Didenko, B.~Silbermann, \emph{Structure of kernels and
cokernels of {T}oeplitz plus {H}ankel operators}, Integral Equations
Oper. Theory, \textbf{10}(1) (2014), 1--31.

\bibitem{DS:2013}
V.\,D. Didenko, B.~Silbermann, \emph{Index calculation for
{T}oeplitz plus {H}ankel operators with piecewise quasi-continuous
generating functions}, Bull. London Math. Soc. \textbf{45}(3)
(2013), 633--650.

\bibitem{Ehr:2004h}
T. Ehrhardt, \emph{Factorization theory for {T}oeplitz+{H}ankel
operators and singular integral operators with flip}, Habilitation
{T}hesis, Technische Universit{\"a}t Chemnitz, 2004.

\bibitem{E2004}
T. Ehrhardt, \emph{Invertibility theory for {T}oeplitz plus {H}ankel
operators and singular integral operators with flip}, J. Funct.
Anal. \textbf{208}(1) (2004), 64--106.

\bibitem{GF1974}
I.\, C. Gohberg, I.\,A. Feldman, \emph{Convolution equations and
projection methods for their solution}, American Mathematical
Society, Providence, R.I., 1974.

\bibitem{RSS:2011}
S. Roch, P.\,A. Santos, B. Silbermann, \emph{Non-commutative
{G}elfand theories. A tool-kit for operator theorists and numerical
analysts}, Universitext, Springer-Verlag London Ltd., London, 2011.

\bibitem{RS1990}
S. Roch, B. Silbermann, \emph{Algebras of convolution operators and
their image in the {C}alkin algebra}, Report MATH, vol.~90, Akademie
der Wissenschaften der DDR Karl-Weierstrass-Institut f\"ur
Mathematik, Berlin, 1990.

\bibitem{RS:2012}
S. Roch, B. Silbermann, \emph{{A handy formula for the Fredholm
index of Toeplitz plus Hankel operators}}, Indag. Math. \textbf{23}
(2012), no.~4, 663--689.

\bibitem{Si:1987}
B. Silbermann, \emph{The {$C^*$}-algebra generated by {T}oeplitz and
{H}ankel operators with piecewise quasicontinuous symbols}, Integral
Equations Operator Theory \textbf{10}(5) (1987), 730--738.



\end{thebibliography}

\end{document}